\documentclass[a4paper,leqno,11pt]{article}
\usepackage{amssymb}
\usepackage{amsmath}
\usepackage[german]{babel}
\begin{document}

\def\refname{ }
\begin{center}
\begin{Large}
\textbf{Algebraic independence of certain Mahler numbers}
\end{Large}
\vskip1.1cm
\large{Keijo V\"a\"an\"anen}
\end{center}
\vskip0.7cm
\noindent
\textit{Abstract}: In this note we prove algebraic independence results for the values of a special class Mahler functions. In particular, the generating functions of Thue-Morse, regular paperfolding and Cantor sequences belong to this class, and we obtain the algebraic independence of the values of these functions at every non-zero algebraic point in the open unit disk. The proof uses results on Mahler's method.\\



\vskip0.4cm
\noindent
\textbf{1. Introduction and results}
\vskip0.2cm
\noindent
In the present paper we are interested in the values of special degree 1 Mahler functions $F(z)$ satisfying a functional equation of the form
\[
p(z) + p_0(z)F(z) + p_1(z)F(z^d) = 0,
\]
where $d \geq 2$ is an integer and $p(z), p_0(z), p_1(z)$ are polynomials satisfying $p_0(z)p_1(z) \neq 0$. The values $F(1/b)$ with integers $b \geq 2$ are called Mahler numbers. The arithmetic properties of such numbers has been an active research area in last years. In a remarkable work Bugeaud \cite{Bu} proved that the irrationality exponent of the Thue-Morse-Mahler numbers $f_{TMM}(1/b)$ equals 2, here 
\[
f_{TMM}(z) = \sum_{n=0}^\infty t_nz^n, 
\]
and $(t_n)$ is the famous Thue-Morse sequence defined recursively by $t_0 = 0, t_{2n} = t_n, t_{2n+1} = 1 - t_n \ (n \geq 0)$. Then similar results were proved  by Coons \cite{C} for the values of the functions
\[
G(z) = \sum_{n=0}^\infty \frac{z^{2^n}}{1 - z^{2^n}}, F(z) = \sum_{n=0}^\infty \frac{z^{2^n}}{1 + z^{2^n}},
\]
and by Guo, Wen and Wu \cite{GWW} and Wen and Wu \cite{WW} for the values of
\[
f_{RPF}(z) = \sum_{n=0}^\infty u_nz^n,\ f_{C}(z) = \sum_{n=0}^\infty v_nz^n, 
\]
respectively, where $(u_n)$ is the regular paperfolding sequence defined by $u_{4n} = 1, u_{4n+2} = 0, u_{2n+1} = u_n \ (n \geq 0)$ and $(v_n)$ is the Cantor sequence on $\{0,1\}$ such that $v_n = 1\ (n \geq 0)$ if and only if the ternary expansion of $n$ does not contain the digit 1. For a unified expression of these (and other) results we refer to \cite{BHWY}. It is also well-known that all these functions obtain transcendental values at every non-zero algebraic point in the open unit disk $\mathbb{D}$. Here our aim is to consider the algebraic independence of these special Mahler numbers.\\

\textbf{Theorem 1}. \textit{For every non-zero algebraic $\alpha \in \mathbb{D}$ the numbers $f_{TMM}(\alpha), f_{RPF}(\alpha)$ and $G(\alpha)$ are algebraically independent over $\mathbb{Q}$. The same holds if we replace $G(\alpha)$ by $F(\alpha)$.}\\

\textbf{Corollary 1}. \textit{For every integer $b \geq 2$, the three numbers
\[
f_{TMM}(\frac{1}{b}),\ f_{RPF}(\frac{1}{b}),\ \sum_{n=0}^\infty \frac{1}{b^{2^n} + 1}
\]
are algebraically independent over $\mathbb{Q}$.}\\

Note that if $b = 2$, then the latest number above is the reciprocal sum of Fermat numbers.

Theorem 1 is obtained from the following more general result. To introduce this, let $d \geq 2$ be a fixed integer, and denote
\begin{equation}\label{1}
T_d(z) = \prod_{n=0}^{\infty}(1 - z^{d^n}),\ U_d(z) = \prod_{n=0}^{\infty}(1 + z^{2d^n}),
\end{equation}
\[
G_{d,j}(z) = \sum_{n=0}^\infty \frac{z^{d^n}}{1 - z^{d^{n+j}}},\ j \in \mathbb{N}_0 := \{0,1,2,\ldots\}.
\]
Then we have\\

\textbf{Theorem 2}. \textit{Let $\alpha \in \mathbb{D}\setminus\!\{0\}$ be an algebraic number. Then the numbers $T_2(\alpha)$ and $G_{2,j}(\alpha), j \in \mathbb{N}_0\setminus\!\{1\}$, are algebraically independent over $\mathbb{Q}$.}\\

\textbf{Theorem 3}. \textit{Let $\alpha \in \mathbb{D}\setminus\!\{0\}$ be an algebraic number. If $d \geq 3$, then the numbers $T_d(\alpha), U_d(\alpha)$ and $G_{d,j}(\alpha), j \in \mathbb{N}_0$, are algebraically independent over $\mathbb{Q}$.}\\

Since
\[
T_2(z) = \frac{1}{1 - z} - 2f_{TMM}(z)
\]
and
\[
G_{2,2}(z) = zf_{RPF}(z),
\]
see \cite{BHWY}, Theorem 2 implies Theorem 1.\\

\textbf{Theorem 4}. \textit{Let $\alpha \in \mathbb{D}\setminus\!\{0\}$ be an algebraic number. Then the numbers
\begin{equation}\label{e1}
T_d(\alpha), U_3(\alpha), G_{3,1}(\alpha), G_{d,j}(\alpha),\ d = 2, 3; j \in \mathbb{N}_0\setminus\!\{1\},
\end{equation}
are algebraically independent over $\mathbb{Q}$.}\\

Since $U_3(z) = f_C(z)$, we immediately obtain the following\\

\textbf{Corollary 2}. \textit{Let $\alpha \in \mathbb{D}\setminus\!\{0\}$ be an algebraic number. Then the numbers $f_{TMM}(\alpha), f_{RPF}(\alpha), F(\alpha)$ and $f_C(\alpha)$ are algebraically independent over $\mathbb{Q}$. In particular, for every integer $b \geq 2$, the four numbers
\[
f_{TMM}(\frac{1}{b}),\ f_{RPF}(\frac{1}{b}),\ \sum_{n=0}^\infty \frac{1}{b^{2^n} + 1}, f_C(\frac{1}{b})
\]
are algebraically independent over $\mathbb{Q}$.}\\ 

To prove Theorems 2, 3 and 4 we consider in Chapter 2 algebraic independence over $\mathbb{C}(z)$ of the functions (\ref{1}). Then Mahler's method can be used to prove our theorems in Chapter 3.

\vskip0.4cm
\noindent
\textbf{2. Algebraic independence of functions}
\vskip0.2cm
\noindent
To study algebraic independence of the functions (\ref{1}) we use the following special case of a result of Kubota \cite{Ku} to be found also in Nishioka \cite [Theorem 3.5]{Ni}.\\

\textbf{Theorem K}. \textit{Let us assume that $f_j(z) \in\mathbb{C}[[z]]\setminus\{0\}\ (j = 0,1,\ldots, m+h)$  converge on $\mathbb{D}$ and satisfy the functional equations
\begin{equation}\label{2}
f_j(z^d) = a(z)f_j(z) + a_j(z),\ j = 0,1,\ldots, m,\ f_{m+i}(z^d) = b_i(z)f_{m+i}(z),\ i = 1,\ldots, h,
\end{equation}
with $a(z), a_j(z), b_i(z) \in \mathbb{C}(z)\setminus\{0\}$. Then the functions $f_j(z)\ (j = 0,1,\ldots, m+h)$ are algebraically independent over $\mathbb{C}(z)$, if $a(z), a_j(z)$ and $b_i(z)$ satisfy the following conditions.}

\noindent
(i) \textit{If $c_0,\ldots, c_m \in \mathbb{C}$ are not all zero, then the functional equation
\[
g(z^d) = a(z)g(z) - \sum_{j=0}^m c_j a_j(z)
\]
does not have a solution $g(z) \in \mathbb{C}(z)$.}

\noindent
(ii) \textit{For any $(n_1,\ldots, n_h) \in\mathbb{Z}^h\setminus\{\underline{0}\}$ the functional equation
\[
r(z^d) = (\prod_{i=1}^{h}b_i(z)^{n_i}) r(z)
\]
does not have a solution $r(z) \in \mathbb{C}(z)\setminus\{0\}$.}\\

The functions (\ref{1}) satisfy functional equations of the form (\ref{2}), namely
\begin{equation}\label{3}
T_d(z^d) = \frac{1}{1 - z}T_d(z),\ U_d(z^d) = \frac{1}{1 + z^2}U_d(z),\ G_{d,j}(z^d) = G_{d,j}(z) - \frac{z}{1 - z^{d^j}},\ j \in \mathbb{N}_0.
\end{equation}
Applying Theorem K we prove first\\

\textbf{Lemma 1}. \textit{If $d \geq 3$, then the functions (\ref{1}) are algebraically independent over $\mathbb{C}(z)$.}\\

\textit{Proof}. Assume, against Lemma 1, that there exists an integer $m \geq 1$ such that the functions $T_d(z), U_d(z)$ and $G_{d,j}(z),\ 0 \leq j \leq m$, are algebraically dependent. We shall prove that these functions satisfy conditions (i) and (ii) of Theorem K and thus obtain a contradiction with this assumption.

Let us consider (i) first. Assume that $c_0, c_1,\ldots, c_m \in \mathbb{C}$ are not all zero. If the functional equation
\begin{equation}\label{4}
g(z^d) = g(z) - \sum_{j=0}^m \frac{c_jz}{1 - z^{d^j}}
\end{equation}  
has a rational solution $g(z)$, then by \cite [Lemma 1]{Ni1} we have
\[
g(z) = \frac{A(z)}{1 - z^{d^m}},
\]
where $A(z)$ is a polynomial. By (\ref{4}), deg $A(z) \leq d^m$. Thus there exist $c \in \mathbb{C}$ and a polynomial $B(z) \neq 0$ with deg $B(z) < d^m$ such that
\[
g(z) = c + \frac{B(z)}{1 - z^{d^m}}.
\]
Letting $z \rightarrow \infty$ in (\ref{4}) we get $c = c + c_0, c_0 = 0$, and
\begin{equation}\label{5}
\frac{B(z^d)}{1 - z^{d^{m+1}}} = \frac{B(z)}{1 - z^{d^m}} - \sum_{j=1}^m \frac{c_jz}{1 - z^{d^j}}.
\end{equation}
We use induction to prove that this is not possible.

If $m = 1$, then (\ref{5}) is of the form
\[
\frac{B(z^d)}{1 - z^{d^2}} = \frac{B(z)}{1 - z^d} - \frac{c_1z}{1 - z^d},\ c_1 \neq 0,
\]
and so
\[
B(z^d) = (B(z) - c_1z)(1 + z^d + z^{2d} +\cdots+ z^{(d-1)d}).
\]
By comparing the coefficients of $z^{kd}$ in this equation we get
\[
B(z) = b_0(1 + z +\cdots+ z^{d-1}) = b_0\frac{1 - z^d}{1 - z}
\]
implying
\[
b_0 = b_0(1 + z +\cdots + z^{d-1}) - c_1z.
\]
Since $d \geq 3$, this leads to a contradiction $b_0 = c_1 = 0$.

Assume now that (\ref{5}) is not possible, if $m$ is replaced by $m - 1 \geq 1$. If we denote
\[
B(z) = \sum_{j=0}^{d^m-1} b_jz^j,
\]
then (\ref{5}) implies
\[
\sum_{j=0}^{d^m-1} b_jz^{dj} = (\sum_{j=0}^{d^m-1} b_jz^j)(1 + z^{d^m} +\cdots+ z^{(d-1)d^m}) - (1 - z^{d^{m+1}})\sum_{j=1}^m \frac{c_jz}{1 - z^{d^j}}.
\]
We compare again the coefficients of $z^{kd}$ in this equation to get $b_0 = b_{d^{m-1}} = b_{2d^{m-1}} =\cdots= b_{(d-1)d^{m-1}}; b_1 = b_{d^{m-1}+1} = b_{2d^{m-1}+1} = \cdots= b_{(d-1)d^{m-1}+1}; \ldots; b_{d^{m-1}-1} = b_{2d^{m-1}-1} =\cdots= b_{d^m-1}$, which means that
\[
B(z) = (\sum_{j=0}^{d^{m-1}-1} b_jz^j)(1 + z^{d^{m-1}} +\cdots+ z^{(d-1)d^{m-1}}) =: B_1(z)\frac{1 - z^{d^m}}{1 - z^{d^{m-1}}}.
\] 
Then, by (\ref{5}),
\[
\frac{B_1(z^d)}{1 - z^{d^m}} = \frac{B_1(z)}{1 - z^{d^{m-1}}} - \sum_{j=1}^m \frac{c_jz}{1 - z^{d^j}}.
\]
If $c_m = 0$, then we have a contradiction by our induction hypothesis. Therefore we necessarily have $c_m \neq 0$, and
\[
B_1(z^d) = B_1(z)(1 + z^{d^{m-1}} +\cdots+ z^{(d-1)d^{m-1}}) - (1 - z^{d^m})\sum_{j=1}^m \frac{c_jz}{1 - z^{d^j}}.
\]
Repeating the above consideration we get
\[
B_1(z) = (\sum_{j=0}^{d^{m-2}-1} b_jz^j)(1 + z^{d^{m-2}} +\cdots+ z^{(d-1)d^{m-2}}) =: B_2(z)\frac{1 - z^{d^{m-1}}}{1 - z^{d^{m-2}}}.
\]
So we have
\[
\frac{B_2(z^d)}{1 - z^{d^{m-1}}} = \frac{B_2(z)}{1 - z^{d^{m-2}}} - \sum_{j=1}^m \frac{c_jz}{1 - z^{d^j}},  
\]
where $c_m \neq 0$. By comparing the poles on both sides of this equation we now get a contradiction.

We next consider the condition (ii). Assume that for some pair $(n_1,n_2) \neq \underline{0}$ the functional equation
\[
r(z^d) = (1 - z)^{-n_1}(1 + z^2)^{-n_2}r(z)
\]
has a rational solution $r(z) \neq 0$, and denote $r(z) = s(z)/t(z)$ with coprime polynomials $s(z)$ and $t(z)$. 

If $n_1, n_2 \geq 0$, then
\[
s(z)t(z^d) = (1 - z)^{n_1}(1 + z^2)^{n_2}s(z^d)t(z).
\]
Since $s(z)$ and $t(z)$ are coprime, this means that $s(z^d)$ is a factor of $s(z)$, and thus $s(z) = s \in \mathbb{C}\setminus\{{0}\}$ and
\[
t(z^d) = (1 - z)^{n_1}(1 + z^2)^{n_2}t(z).
\]
Since the polynomials $t(z)$ and $t(z^d)$ have the same multiplicity of zero at $z = 1$, we necessarily have $n_1 = 0$ and $(d-1)D = 2n_2$, where $D := $deg $t(z)$. If $d \geq 4$, then $D < n_2$, and so the equation
\begin{equation}\label{e2}
t(z^d) = (1 + z^2)^{n_2}t(z)
\end{equation}
is not possible. If $d = 3$, then $D = n_2$. The equation $z^3 = c \in \mathbb{C}$ may have at most one of $i$ or $-i$ as a root. From this it follows that (\ref{e2}) is not possible, if $d = 3$. The case $n_1, n_2 \leq 0$ is similar.

If $n_1, -n_2  \geq 0$, we denote $N := -n_2$, and then
\[
s(z^d)t(z)(1 - z)^{n_1} = (1 + z^2)^N s(z)t(z^d).
\]
Thus there exists a polynomial $u(z)$ such that $t(z^d)u(z) = t(z)(1-z)^{n_1}$. Since $t(z)$ and $t(z^d)$ have the same multiplicity of zero at $z = 1$, we obtain  $u(z) = (1-z)^{n_1}v(z)$ with some polynomial $v(z) \neq 0$. But then $t(z^d)v(z) = t(z)$ giving $v(z) = 1, t(z) = t \in \mathbb{C}\setminus\{0\}$. So $s(z^d)(1-z)^{n_1} = (1+z^2)^N s(z)$, and using again the fact that $s(z)$ and $s(z^d)$ have the same multiplicity of zero at $z = 1$, we necessarily have $n_1 = 0$. But then $s(z^d) = (1+z^2)^Ns(z)$, and this is analogous to (\ref{e2}) and so impossible, as we saw above. The case $-n_1, n_2 \geq 0$ can be considered in a similar way. This proves (ii).

Theorem K gives now the truth of Lemma 1.\\

\textbf{Lemma 2}. \textit{The functions $T_2(z), G_{2,j}(z), j \in \mathbb{N}_0\setminus\!\{1\}$, are algebraically independent over $\mathbb{C}(z)$.}\\

\textit{Proof}. The proof of the induction step in condition (i) above works also in the case $d = 2$, but the starting point of the induction does not hold, since $G_{2,1}(z) = z/(1-z)$ is a rational function, see \cite [Theorem 9]{DN}. So in this case we delete $G_{2,1}(z)$, start with $m = 2$ and consider the functional equation
\[
\frac{B(z^2)}{1 - z^8} = \frac{B(z)}{1 - z^4} - \frac{c_2z}{1 - z^4},\ c_2 \neq 0.
\]  
Let $B(z) = b_0 + b_1z + b_2z^2 + b_3z^3$. Then
\[
b_0 + b_1z^2 + b_2z^4 + b_3z^6 = (b_0 + b_1z + b_2z^2 + b_3z^3)(1 + z^4) +c_2z(1 + z^4),
\]
and by comparing the coefficients of even powers of $z$ on both sides we have $B(z) = (b_0 + b_1z)(1 + z^2)$. Thus
\[
\frac{b_0 + b_1z^2}{1 - z^4} = \frac{b_0 + b_1z}{1 - z^2} - \frac{c_2z}{1 - z^4},
\]
and so $b_0 + b_1z^2 = (b_0 + b_1z)(1 + z^2) - c_2z$ implying $b_0 = b_1 = c_2 = 0$, a contradiction. So we may now start the induction from $m = 2$ and continue as in the proof of Lemma 1 to obtain the condition (i).

To consider the condition (ii), let us assume that for some integer $n \neq 0$ the functional equation
\[
r(z^2) = (1 - z)^{-n}r(z)
\]
has a rational solution $r(z) \neq 0$, and denote $r(z) = s(z)/t(z)$ with coprime polynomials $s(z)$ and $t(z)$. If $n > 0$, then
\[
s(z^2)(1-z)^nt(z) = s(z)t(z^2).
\]
Since $s(z)$ and $t(z)$ are coprime, this means that $s(z^2)$ is a factor of $s(z)$, and thus $s(z) = s \in \mathbb{C}\setminus\{{0}\}$. Therefore $(1-z)^nt(z) = t(z^2)$, which leads immediately to a contradiction $n = 0$. The case $n < 0$ is similar.

Thus Lemma 2 is true.

We note that we cannot include $U_2(z)$ to the functions in Lemma 2, since $U_2(z) = 1/(1-z^2)$.

\vskip0.4cm
\noindent
\textbf{3. Proof of Theorems 2, 3 and 4}
\vskip0.2cm
\noindent
 
We shall need the following basic result of Mahler's method given in \cite [Theorem 4.2.1]{Ni}.\\

\textbf{Theorem N1}. \textit{Let $K$ denote an algebraic number field. Suppose that $f_1(z),...,f_m(z)\in K[[z]]$ converge
in some disk $U\subset\mathbb{D}$ about the origin, where they satisfy the matrix functional equation
\[
{}^\tau(f_1(z^d),...,f_m(z^d))=\mathcal{A}(z)\cdot{}^\tau(f_1(z),...,f_m(z))+
{}^\tau(b_1(z),...,b_m(z))
\]
with $\mathcal{A}(z)\in\mathrm{Mat}_{m\times m}(K(z)), \tau$ indicating the matrix transpose, and $b_1(z),...,b_m(z)\in K(z)$. If
$\alpha\in U\setminus\{{0}\}$ is an algebraic number such that none of the $\alpha^{d^j} \,(j\in\mathbb{N}_0)$ is a pole of $b_1(z),...,b_m(z)$ and the entries of $\mathcal{A}(z)$, then the following inequality holds}
\[
\mathrm{trdeg}_\mathbb{Q}\mathbb{Q}(f_1(\alpha),...,f_m(\alpha))\ge
\mathrm{trdeg}_{K(z)}K(z)(f_1(z),...,f_m(z)).\\
\]

This result with Lemmas 1 and 2 gives immediately the truth of Theorems 2 and 3.

In Theorem 4 we consider two different values $d = 2$ and $d = 3$. For this we recall the following special case of \cite [Theorem 1]{Ni2}, where $f_{i,1}(z),\ldots, f_{i,m_i}(z) \in K[[z]]\ (i = 1,2)$ converge in $\mathbb{D}$ and satisfy
\[
f_{i,j}(z) = a_{i,j}(z)f_{i,j}(z^{d_i}) + b_{i,j}(z),\ i = 1, 2; j = 1,\ldots, m_i,
\]
with $a_{i,j}(z), b_{i,j}(z) \in K(z)$ and $a_{i,j}(0) = 1$.\\

\textbf{Theorem N2}. \textit{Suppose that $\log d_1/\log d_2 \notin \mathbb{Q}$. Let $\alpha\in \mathbb{D}\setminus\{{0}\}$ be an algebraic number such that none of the $\alpha^{d_i^k} \,(i = 1,2; k\in\mathbb{N}_0)$ is a pole of $a_{i,j}(z), b_{i,j}(z)$ and $a_{i,j}(\alpha^{d_i^k}) \neq 0$. If, for both values $i = 1,2$, the functions $f_{i,1}(z),\ldots, f_{i,m_i}(z)$ are algebraically independent over $K(z)$, then the values
\[
f_{i,j}(\alpha),\ i = 1,2; j = 1,\ldots, m_i,
\]
are algebraically independent over $\mathbb{Q}$.}\\

We now choose $K = \mathbb{Q}, d_1 = 2, d_2 = 3$, and assume, against Theorem 4, that the numbers (\ref{e1}) are algebraically dependent. Then there exists an integer $m \geq 2$ such that the numbers
\begin{equation}\label{6}
T_d(\alpha), U_3(\alpha), G_{d,0}(\alpha), G_{3,1}(\alpha), G_{d,j}(\alpha),\ d = 2,3; j = 2,\ldots,m,
\end{equation}
are algebraically dependent. By Lemma 2 the functions $T_2(z), G_{2,0}(z), G_{2,j}(z)\ (j = 2,\ldots,m)$ are algebraically independent over $K(z)$, and, by Lemma 1, also the functions $T_3(z), U_3(z),\\ G_{3,j}(z)\ (j = 0,1,\ldots,m)$ are algebraically independent over $K(z)$. Thus Theorem N2 implies the algebraic independence of the numbers (\ref{6}). This contradiction proves the truth of Theorem 4.

\vskip5cm
\noindent
\textbf{References}
\vskip0.2cm
\vspace{-1,5cm}

\vskip0.8cm
\begin{small}
\noindent
Author's address:   \\
\vskip0.1cm
\noindent
Keijo V\"a\"an\"anen  \\
Department of Mathematical Sciences \\ 
University of Oulu \\
P. O. Box 3000  \\
90014 Oulun Yliopisto, Finland \\
E-mail: keijo.vaananen@.oulu.fi \\
\end{small}


\begin{thebibliography}{9}

\bibitem{Bu} Yann Bugeaud, On the rational approximation to the Thue{-}Morse{-}Mahler numbers,
\textit{Ann. Inst. Fourier (Grenoble)}, \textbf{61}(5) (2011) 2065--2076.

\bibitem{BHWY} Yann {Bugeaud}, Guo-Niu {Han}, Zhi-Ying {Wen}, and Jia-Yan {Yao}, Hankel determinants, Pad\'{e} approximations, and irrationality exponents, arXiv:1503.02797 (2015).

\bibitem{C} M. Coons, On the rational approximation of the sum of the reciprocals of the Fermat numbers, \textit{The Ramanujan Journal} \textbf{30}, no. 1 (2013) 39-65.

\bibitem{DN} D. Duverney and Ku. Nishioka, An inductive method for proving transcendence of certain series,
\textit{Acta Arith.} \textbf{110} (2003) 305--330.

\bibitem{GWW} Ying-Jun {Guo}, Zhi-Xiong {Wen}, and Wen {Wu}, On the irrationality exponent of the regular paperfolding numbers, \textit
{Linear Algebra Appl.}, \textbf{446} (2014) 237 -- 264.

\bibitem{Ku} K. K. Kubota, On the algebraic independence of holomorphic solutions of
certain functional equations and their values, \textit{Math. Ann.} \textbf{227} (1977) 9--50.

\bibitem{Ni2} Ku. Nishioka, Algebraic independence by Mahler's method and S-unit equations, \textit{Compositio Math.}, \textbf{92} (1994) 87--110.

\bibitem{Ni} Ku. Nishioka, \textit{Mahler Functions and Transcendence}, LNM 1631 (Springer, Berlin et al., 1996).

\bibitem{Ni1} Ku. Nishioka, Algebraic independence of reciprocal sums of binary recurrences, \textit{Monatsh. Math.} \textbf{123} (1997) 135--148.

\bibitem{WW} Zhi-Xiong {Wen} and Wen {Wu}, Hankel determinants of the Cantor sequence, \textit{Scientia Sinica Mathematica (Chinese)}, \textbf{44} (10) (2014) 1059--1072.

\end{thebibliography}
\end{document}